\documentclass{amsart}
\usepackage{amsmath, amssymb, amsfonts, url, color, hyperref}
\usepackage[foot]{amsaddr}
\newtheorem{theorem}{Theorem}[section]
\newtheorem{lemma}{Lemma}[section]

\newtheorem{question}{Question}[section]
\newtheorem*{question*}{Question}
\newtheorem*{problem*}{Problem}
\newtheorem{definition}{Definition}[section]

\newtheorem*{acknowledgements*}{Acknowledgements}
\numberwithin{equation}{section}

\newcommand{\CC}{\mathbb{C}}
\newcommand{\RR}{\mathbb{R}}

\begin{document}

\title[A survey on generalizations of Forelli's theorem]{A survey on generalizations of Forelli's theorem and related pluripotential methods}

\author{Ye-Won Luke Cho}

\begin{abstract}
	We present a survey on recent developments of generalizations of Forelli's analyticity theorem and related pluripotential methods. 
\end{abstract}

\maketitle

\section{Introduction}

When the higher-dimensional complex analysis began to emerge in the early 1900s, a natural concern was to  establish criteria for the complex analyticity of functions of several complex variables. Then the study of analyticity theorems started with the following celebrated 

\begin{theorem}[Hartogs \cite{Hartogs1906}]\label{Hartogs-original}
 A complex-valued function $f:\Omega\to \mathbb{C}$  on an open set $\Omega\subset \CC^n$ is holomorphic if $f$ is separately holomorphic, i.e., it satisfies the Cauchy-Riemann equations in each complex variable separately.
\end{theorem}

 Note that the converse of the theorem is trivial. Earlier, Osgood \cite{Osgood} proved Theorem \ref{Hartogs-original} under an additional assumption that the given function is continuous. The ingenious idea of Hartogs was to use subharmonic functions (Lemma \ref{Hartogslemma}) to show that such an assumption is redundant. Since then, the theory of plurisubharmonic functions which we now call the \textit{pluripotential theory} has been developed by various researchers. For the applications of the modern pluripotential methods to generalizations of Theorem \ref{Hartogs-original}, see \cite{JP11} and references therein.
 
The following theorem of Forelli should be second only to Hartogs' analyticity theorem among various complex analyticity theorems.

\begin{theorem}[Forelli \cite{Forelli77}] \label{Forelli-original}
	If a function $f\colon B^n \to \CC$ defined on the open unit ball $B^n\subset \mathbb{C}^n$ satisfies the following 
	two conditions
	\begin{enumerate}
		\setlength\itemsep{0.3em}
		\item $f\in C^{\infty}(0)$, meaning that for any positive integer $k$ there
		exists an open neighborhood $V_k$ of the origin $0$ such that $f \in C^k 
		(V_k)$, and
		\item the correspondence $f_v: z\in B^1\to f(z v)$ is holomorphic on $B^1$ for each
		$v \in \CC^n$ with $\|v\|=1$,
	\end{enumerate}
	then $f$ is holomorphic on $B^n$.
\end{theorem}

   At first glance, one may think that Forelli's theorem can be obtained as an easy corollary of Hartogs' theorem via the birational blow-up of $B^n$ at the origin. Although the idea can be realized as shown in \cite{Kim13}, the proof requires additional nontrivial analysis. 
   
   Naturally, there have been many attempts to generalize Theorem \ref{Forelli-original}. But it turned out that Condition (1) cannot be weakened to finite differentiability as numerous counterexamples were found; see \cite{JKS16}. Recently, Condition (2) has been generalized successfully to various directions, starting with the work of Chirka \cite{Chirka06}. The purpose of this article is to present the recent developments of generalizations of Forelli's analyticity theorem and related pluripotential methods. The contents of the article are based upon the author's talk at the POSTECH Conference on Complex Analytic Geometry held in Pohang, South Korea for the period 18-22 July 2022.

\subsection*{Acknowledgement}
 Most parts of the article were written while the author was supported by the National Research Foundation of Korea (NRF-2018R1C1B3005963, NRF-2021R1A4A1032418). The author is currently supported by the National Research Foundation of Korea (NRF-2021R1A4A1032418, NRF-
 2023R1A2C1007227, RS-2023-00246259).

\section{The original proof of Forelli}
The original version of Forelli's theorem is concerned with the pluriharmonicity of functions harmonic along linear complex discs passing through the origin. But with minor adjustments, Forelli's proof \cite{Forelli77} also yields Theorem \ref{Forelli-original} as noted in \cite{Stoll80}. We shall recapitulate the proof in this section.

\medskip
\textit{Sketch of the proof of Theorem \ref{Forelli-original}}. Let $f$ be the given function. Note that, by Condition (2), one can define 
\[
f_m(z):=\frac{1}{2\pi}\int_{0}^{2\pi}f(ze^{i\theta})\,e^{-im\theta}d\theta
\]
for each nonnegative integer $m\geq 0$ and $z\in B^n$. Let $S_f:=\sum_{m=0}^{\infty}f_m$ be the formal sum of the sequence $\{f_m\}$. Fix $z\in B^n$ and choose $t\in \mathbb{C}$ with $|t|\cdot \|z\|<1$. As the correspondence $t \mapsto f(tz)$ is holomorphic by the assumption, we have
\[
f(tz)=\sum_{m=0}^{\infty}c_m(z)\cdot t^m
\]
for some sequence $\{c_m(z)\}\subset \mathbb{C}$. Then it follows from the definition of $f_m$ that $f_m(tz)=t^mc_m(z)$. Letting $t \to 1$, we obtain $f_m(z)=c_m(z)$. So $f\equiv S_f$ on $B^n$ and
\begin{equation}\label{homogeneous}
	f_m(tz)=t^mf_m(z)~\text{for each}~t\in B^1,\,z\in B^n.
\end{equation}
Since $f\in C^{\infty}(0)$, we may differentiate the terms in (\ref{homogeneous}) $m$ times with respect to the complex variable $t$. Then one easily verifies that each $f_m$ is a holomorphic  homogeneous polynomial of degree $m$. Hence, we say that the formal Taylor series $S_f$ is of $\textit{holomorphic type}$. Note also that $f$ is continuous on a neighborhood of the origin. So there are numbers $2r\in (0,1),\, M>0$ such that $|f(z)|\leq M$ if  $z\in B^n(0;2r) := \{z \in \CC^n \colon \|z-a\|<2r \}$. Then for each $z\in B^n(0;r)$, we have
\begin{align*}
|f_m(z)|&=\frac{1}{2^m}\cdot |f_m(2z)|\leq \frac{1}{2\pi}\frac{1}{2^m}\cdot \int_{0}^{2\pi}|f(2ze^{i\theta})|\,d\theta\leq \frac{M}{2^m}.
\end{align*}
Therefore, $f\equiv S_f$ is holomorphic on $B^n(0;r)$ by the Weierstrass $M$-test. To prove that $f$ is holomorphic on $B^n$, we recall the following celebrated lemma of Hartogs.
\begin{lemma}[Hartogs \cite{Hartogs1906}]\label{Hartogslemma}
	Let $\{u_m\}$ be a sequence of subharmonic functions on an open set $\Omega\subset \mathbb{C}^n$ and $C\in \mathbb{R}$ a constant such that
	\begin{enumerate}
		\setlength\itemsep{0.1em}
		\item $\{u_m\}$ is locally uniformly bounded from above on $\Omega$, and
		\item $\limsup\limits_{m\to \infty}u_m(z)\leq C$ for any $z\in \Omega$.
	\end{enumerate}
	If $K$ is a compact subset of $\Omega$ and $\epsilon$ is a positive number, then there exists a positive integer $N=N(K,\epsilon)$ such that $u_m(z)\leq C+\epsilon$ whenever $m\geq N$ and $z\in K$.
\end{lemma}
\noindent
Let $\{u_m\}$ be a sequence of subharmonic functions on $\mathbb{C}^n$ defined as $u_m(z):=|f_m(z)|^{\frac{1}{m}}$. Then applying Lemma \ref{Hartogslemma} to the sequence $\{u_m\}$, one can conclude that $S_f$ converges uniformly on each compact subset of $B^n$; see p.362 of \cite{Forelli77}. \hfill $\Box$

\section{Forelli's theorem at the level of formal power series}\label{Section3}

Recall that the proof of Theorem \ref{Forelli-original} was carried out in the following two steps:

 \smallskip
\begin{narrower}
	
	\textbf{Step 1.} The formal Taylor series $S_f$ of $f$ is of holomorphic type.
	
	\textbf{Step 2.} The formal series $S_f$ converges uniformly on some $B^n(0;r),~r<1$.
	
\end{narrower}
\noindent
 Then by Lemma \ref{Hartogslemma},  $f\equiv S_f$ is holomorphic on $B^n$. In this section, we shall introduce several studies regarding \textbf{Step 2}.
\subsection{Solutions of Bochner's problem by Zorn, Ree, Lelong, and Cho-Kim}
A line of research concerning the convergence of formal power series of holomorphic type originates from the following question of Bochner:

  \begin{question}[Bochner]
	Let $S=\sum a_{ij}z_1^iz_2^j$ be a formal power series with complex coefficients such that every substitution of convergent power series with 
	complex coefficients $z_1=\sum b_it^i$, $z_2=\sum c_it^i$ produces a 
	convergent power series in $t$.  Is $S$ convergent on some neighborhood 
	of $0\in \mathbb{C}^2?$
  \end{question}

The question was answered affirmatively by Zorn \cite{Zorn47}. He proved that $S\in \mathbb{C}[[z_1,z_2]]$ is uniformly convergent on an open neighborhood of the origin if the map $t\in \mathbb{C}\to S(at,bt)$ has a nonvanishing radius of convergence $R=R(a,b)>0$ for each $(a,b)\in \mathbb{C}^2$. Zorn also remarked in \cite{Zorn47} that it would be interesting to know whether his result can be generalized to the real case. Ree \cite{Ree49} clarified the meaning of the `real case' and obtained the same conclusion when the set of linear complex discs in Zorn's theorem is replaced by the set $\{t\in \mathbb{C} \to (at,bt):(a,b)\in \mathbb{R}^2\}$. Motivated by the works of Zorn and Ree, Lelong \cite{Lelong51} introduced the following
\begin{definition}[\cite{Lelong51}]\label{definitionLelong}
\normalfont
	A set $E\subset \mathbb{C}^2$ is called \textit{normal} if any formal power series $S\in \mathbb{C}[[z_1,z_2]]$ enjoying the property that $S_{a,b}(t):=S(at,bt)\in \CC[[t]]$ has a positive radius of convergence 	for every $(a,b)\in E$ becomes holomorphic on some open neighborhood of the origin in $\mathbb{C}^2$.
\end{definition}

In the terminology of Lelong, what the theorems of Zorn and Ree say is that $\mathbb{R}^2$ and $\mathbb{C}^2$ are normal sets, respectively. In \cite{Lelong51}, Lelong also characterized a normal set in $\mathbb{C}^2$ using the potential theory on the complex plane. This was generalized to a higher-dimensional principle by Cho-Kim in \cite{CK21}. To state the analyticity theorem of Cho-Kim, recall that a set $F\subset \mathbb{C}^n$ is $\textit{pluripolar}$ if there is a nonconstant plurisubharmonic function $u$ on $\mathbb{C}^n$ such that $u\equiv -\infty$ on $F$. For each nonempty set $F\subset \mathbb{C}^n$, the $\textit{direction set}$ of $F$ is defined as
\[
F':=\bigg\{\bigg(\frac{z_2}{z_1},\dots, \frac{z_n}{z_1}\bigg)\in \mathbb{C}^{n-1}:(z_1,\dots,z_n)\in F,~z_1\neq 0\bigg\}.
\] 
Obviously, one can extend Definition \ref{definitionLelong} to higher dimensions. Then now we can state the following

\begin{theorem}[Cho-Kim \cite{CK21}]\label{CK21}
A set $F\subset \mathbb{C}^n$ is normal if $F'\subset \mathbb{C}^{n-1}$ is not pluripolar.  
\end{theorem}

Here, we give a sketch of the proof of the theorem when $n=2$. Let 
\[
S(z_1,z_2)= \sum a_{i,j}{z_1}^{i}{z_2}^{j}
\]
be a formal power series for which $S_{a_1,a_2}(t):=S(a_1t,a_2t)$ has a positive radius of 
convergence $R_{(a_1,a_2)}>0$ for every $(a_1,a_2)\in F$. Then we are to show that $S$ is holomorphic on some open neighborhood 
of $0$ in $\mathbb{C}^2$.  Note that, for any $b\in F'$, $S_{1,b}$ converges absolutely and uniformly on $\frac{1}{2}R_{(1,b)}$. So it can be rearranged as follows:
\[
S_{1,b}(t)=S(t,bt)=\sum_{m=0}^{\infty}\bigg(\sum_{j=0}^{m}a_{m-j,j}b^j\bigg)t^m = \sum_{m=0}^{\infty}P_m(b)\cdot t^m,
\]
where
\[
P_m(z):=\sum_{j=0}^{m}a_{m-j,j}z^j\in \mathbb{C}[z],~\textup{deg}\,P_m\leq m.
\]
 Then the root test implies the following pointwise estimate:
\begin{equation}\label{limsup}
\limsup\limits_{m\to \infty}\frac{1}{m}\,\textup{log}\,|P_m(b)|<\infty\,\text{for each} ~b\in F'.
\end{equation}
Note that each function $\frac{1}{m}\,\textup{log}\,|P_m(z)|$ is subharmonic on $\mathbb{C}.$ The crux of the proof turns out to be the following version of Lemma \ref{Hartogslemma}; see Proposition 4.1 in \cite{CK21}.

\begin{theorem}[Lelong \cite{Lelong51}, Cho-Kim \cite{CK21}]\label{estimate}
Let $\{P_m\}\subset \mathbb{C}[z_1,\ldots,z_n]$ be a sequence of 
polynomials with $\textup{deg}\,P_m\leq m$ for each positive integer $m$. If $F\subset \mathbb{C}^n$ is nonpluripolar and
\[
\limsup\limits_{m\to \infty}\frac{1}{m}\, \textup{log}\,|P_m(z)|<\infty\,~\text{for each}~ z\in F,
\]
then for each compact subset $K$ of $\mathbb{C}^n$, there exists a constant $M=M(K)>0$ such that $\frac{1}{m}\,\textup{log}\,|P_m(z)|<\textup{log}\,M$ for any $z\in K$, $m\geq 1$.
\end{theorem}
\noindent
So we obtain from (\ref{limsup}) that $|P_m(z)|\leq M^m$ for each $m\geq 1$ and $z\in K:=\bar{B}^n$. Then the conclusion follows from the Cauchy estimate and the Weierstrass $M$-test.

 We remark that the proof of Theorem \ref{estimate} had to await the celebrated theorem of Bedford-Taylor on the equivalence of pluripolar sets and negligible sets in $\mathbb{C}^n$; see Theorem 7.1 in \cite{BedfordandTaylor82}. If $n=1,$ then the theorem of Bedford-Taylor reduces to the result of H. Cartan \cite{Cartan42} used in \cite{Lelong51} (cf. $[a_3]$ on p.14 of \cite{Lelong51}).
\subsection{Projective capacity and related extremal plurisubharmonic function}

The classical theorem of Abel says that, if a formal series $S\in \mathbb{C}[[z]]$ converges at $z_0\in \mathbb{C}-\{0\}$, then $S$ converges uniformly on $B^1(0;r)$ for any $0<r<|z_0|$. However, this does not hold for formal series of several complex variables; Leja showed that there exists a formal power series $S\in \mathbb{C}[[z_1,z_2]]$ that converges pointwise on a given countable set $E\subset \mathbb{C}^2$ but does not converge uniformly on any open neighborhood of the origin. Then it would be natural to address the following
\begin{problem*}[Leja]
Characterize a subset $E$ of $\mathbb{C}^n$ for which the following statement holds: any formal series $S\in \mathbb{C}[[z_1,\dots,z_n]]$ that  converges at each point of $E$ must converge uniformly on an open neighborhood of the origin.
\end{problem*}

We refer the reader to \cite{Siciak90} for a detailed historical account of the problem of Leja. In \cite{Siciak90}, Siciak gave a solution to the problem of Leja using the notion of projective capacity and related plurisubharmonic extremal function. Then the solution also yields a complete characterization of $F_{\sigma}$ normal sets as noted by Levenberg and Molzon \cite{LevenMol88}.  To recapitulate the solution of Siciak, we first recall the following

\begin{definition}[\cite{Siciak82}]\label{defofext}
	\normalfont
    The set of $\textit{homogeneous plurisubharmonic functions}$ on $\mathbb{C}^n$ is defined as
	\[
	H:=\{u\in \textup{PSH}(\mathbb{C}^n): u\geq 0~\text{on}~\mathbb{C}^n, \,u(tz)=|t|u(z)\; \forall z\in \mathbb{C}^n,~ t\in \mathbb{C}\}.
	\]
	For each bounded subset $E$ of $\mathbb{C}^n$, define
	\[
	\Psi_E(z):=\textup{sup}\,\{u(z):u\in H,\; u\leq 1~\text{on}~E\},\;~\forall z\in \mathbb{C}^n.
	\]
	If $E$ is unbounded, then we set
	\[
	\Psi_{E}(z):=\textup{inf}\,\{\Psi_{F}(z): F\subset E~\text{is bounded} \},\;~\forall z\in \mathbb{C}^n.
	\]	
	Denote by $S^{2n-1} := \{z \in \CC^{n}=\RR^{2n} \colon \|z\|=1\}.$ 	The $\textit{projective capacity}$ of a set $E\subset \mathbb{C}^n$ is defined as
	\[
	\rho(E):=\textup{inf}\,\{\|u\|_{E}:u\in H,\, \|u\|_{S^{2n-1}}=1\},~\text{where}~\|u\|_{E}:=\sup_{z\in E}|u(z)|.
	\]
\end{definition}

 Recall that a set $E\subset \mathbb{C}^n$ is $\textit{circular}$ if $(e^{i\theta}z_1,\dots,e^{i\theta}z_n)\in E$ for any $\theta \in \mathbb{R}$ and $(z_1,\dots,z_n)\in E$. From the potential theoretic point of view, the importance of the family $H$ comes from the fact that a circular set $E\subset \mathbb{C}^n$ is pluripolar if, and only if, there is a function $u\in H$ such that $u\equiv 0$ on $E$. If $u\in H$ and $u\equiv 0$ on $E$, then $\rho(E)=0$ by the definition and $\Psi^{\ast}_{E}\equiv +\infty$ as $m\cdot u\leq \Psi_{E}$ on $\mathbb{C}^n$ for each $m>0$. Here, the asterisk denotes the upper-semicontinuous regularization of the function. In the construction of $u\in H$ vanishing on a given circular pluripolar set, Lemma \ref{Hartogslemma} again plays an important role; see Proposition 2.20 in \cite{Siciak82}. If a circular set $E$ is nonpluripolar, then it is known that $\rho(E)\neq 0$ and $\Psi^{\ast}_{E}\in H$.
 
 For each nonempty set $F\subset S^{2n-1}$, define a circular set $S_0(F):=\{zv:z\in B^1,\,v\in F\}$. Let $\psi(z,v):=zv\in S_0(F)$ for each $z\in B^1$ and $v\in F$. The pair $(S_0(F),\psi)$ will be referred as a $\textit{suspension of linear discs}$ at the origin. For simplicity, we will identify a suspension with its underlying set. One can check that $S_0(F)$ is nonpluripolar if, and only if, $F'$ is nonpluripolar. Then now we can state the analyticity theorem of Siciak in the following form:

\begin{theorem}[Siciak \cite{Siciak90}]\label{Siciak original}
Let $F\subset S^{2n-1}$ be a set such that $F'\subset \mathbb{C}^{n-1}$ is nonpluripolar. If $S\in $ $\mathbb{C}[[z_1,\dots,z_n]]$ is a formal series for which the map $z\in B^1\to S(zv)$ is holomorphic for any $v\in F,$ then $S$ is holomorphic on a domain of holomorphy
\[
\Omega:=\big\{z\in \mathbb{C}^n:\Psi^{\ast}_{S_0(F)}(z)<1\big\}\supset B^n(0;\rho(S_0(F)))
\]
containing the origin. Conversely, if $F\subset S^{2n-1}$ is a circular $F_{\sigma}$ set such that $F'$ is pluripolar, then there exists a formal series $S\in \mathbb{C}[[z_1,\dots,z_n]]$ such that the correspondence $z\in B^1 \to S(zv)$ is holomorphic for each $v\in F$ but $S$ does not converge uniformly on any open neighborhood of the origin.
\end{theorem} 

  To prove that the formal series $S$ is holomorphic on $\Omega$, we first write it as a formal sum $S=\sum_{m=0}^{\infty}q_m$ of holomorphic homogeneous polynomials $\{q_m\}$ with $\textup{deg}\,q_m=m$. The crux of the proof is the following $\textit{Bernstein-Walsh type inequality}$ which follows from Definition \ref{defofext}:
\begin{equation}\label{BWinequality}
|q_m(z)|\leq \|q_m\|_{E}\cdot \{\Psi^{\ast}_{E}(z)\}^m~\text{for each}~z\in \mathbb{C}^n.
\end{equation}
It should be noted that the inequality becomes trivial if $E\subset \mathbb{C}^n$ is a pluripolar circular set by the preceding arguments. Then using the Cauchy estimate of $S$ along each complex line in the suspension, one can obtain the following estimate from (\ref{BWinequality}):
\[
\limsup_{m\to \infty}|q_m(z)|^{\frac{1}{m}}\leq \Psi^{\ast}_{S_0(F)}(z)~\text{for any}~z\in \mathbb{C}^n.
\]
Therefore, $S$ is holomorphic on $\Omega$ by the root test. We refer the reader to the proof of Theorem 3.1 in \cite{Siciak90} for the details. We also remark that the region of convergence of $S$ can be estimated using the pluricomplex Green function with pole at infinity as in \cite{Sadullaev22}. 
\section{Generalizations of Forelli's analyticity theorem}
Recent generalizations of Forelli's analyticity theorem started with the work of Chirka \cite{Chirka06}. He proved that the radial complex discs in Theorem \ref{Forelli-original} can be replaced with a singular foliation of $B^2$ by holomorphic curves transversal at the origin to obtain the same conclusion. This was generalized to a higher-dimensional principle by Joo-Kim-Schmalz in \cite{JKS13} and finally by Cho-Kim \cite{CK21} to the case where the foliation is parametrized by an open subset of $S^{2n-1}$.

In another direction, Kim-Polestky-Schmalz observed in \cite{KPS09} that each radial complex disc passing through the origin in $\mathbb{C}^n$ is an integral curve of the complex Euler vector field $E=\sum_{k=1}^{n}z_k\frac{\partial}{\partial z_k}$. Then they proved that the set of integral curves of a \textit{diagonalizable vector field} of the form
\[
X=\sum_{k=1}^{n}\lambda_kz_k\frac{\partial}{\partial z_k},~\{\lambda_1,\dots,\lambda_n\}\subset \mathbb{C}
\]
can replace the set of complex lines in Theorem \ref{Forelli-original} if, and only if, the field is \textit{aligned}, i.e., $\lambda_j/\lambda_k>0$ for each $j,k\in \{1,\dots,n\}$. This was generalized to the case of nondiagonalizable vector fields contracting at the origin in \cite{JKS16}. We say that a holomorphic vector field $X$ defined on an open neighborhood of the origin in $\mathbb{C}^n$ is  $\textit{contracting at the origin}$ if the flow-diffeomorphism $\Phi_t$ of $\textup{Re}\,X$ for some $t<0$ satisfies: (1) $\Phi_t(0)=0$, and (2) every eigenvalue of the matrix $d\Phi_t|_{0}$ has absolute value less than 1. By the Poincar\'{e}-Dulac theorem, there exists a local holomorphic coordinate system near the origin such that $X$ takes the following form:
\begin{equation}\label{contractingVF}
	X=\sum_{k=1}^{n}(\lambda_kz_k+g_k(z))\frac{\partial}{\partial z_k},
\end{equation}
where $g_k\in \mathbb{C}[z_1,\dots,z_n]$ and $\lambda_k\in \mathbb{C}$ for each $k$. $X$ is said to be \textit{aligned} if $\lambda_j/\lambda_k>0$ for each $j,k\in \{1,\dots,n\}$. By replacing the complex time $t$ with $\zeta:=\lambda_1\cdot t$ in the complex flow map $\Phi^X(z,t)$, we will assume that $\lambda_i>0$ for each $i$ whenever $X$ is aligned.

Note that the aforementioned two theories complement each other in a sense; the partial foliation considered by Cho-Kim need not be generated by a vector field, whereas the integral curves of a holomorphic vector field need not intersect mutually transversally at origin. In this section, we shall briefly examine the ideas in the proofs of the generalizations of Theorem \ref{Forelli-original}.
\subsection{Foliation viewpoint}
To formulate a family of holomorphic curves intersecting mutually transversally at a point $p$ of a domain $\Omega\subset \mathbb{C}^n$, Cho-Kim introduced the notion of \textit{$C^1$ pencil of holomorphic discs} at $p\in\Omega$ in \cite{CK21}. The primary example of such a pencil is the \textit{standard pencil} $S_0(U)$ of linear discs at the origin in $\mathbb{C}^n$ parametrized by an open subset $U$ in $S^{2n-1}$. Then a general pencil is defined as a $C^1$ diffeomorphic image of a standard pencil into $\Omega$. Here, we require that the diffeomorphism (1) sends the origin to the point $p$, and (2) deforms each linear disc holomorphically to a Riemann surface in $\Omega$.

Let $\Omega\subset \mathbb{C}^n$ be a domain admitting a $C^1$ pencil of holomorphic discs at $p\in \Omega$. Then the theorem of Cho-Kim \cite{CK21} says that, if a function $f:\Omega\to \mathbb{C}$ is (1) smooth at $p$, and (2) holomorphic along each Riemann surface of the given pencil, then $f$ is holomorphic on the union of an open ball centered at $p$ and the underlying (open) set of the pencil. If $U=S^{2n-1}$, then the theorem reduces to the result of Joo-Kim-Schmalz in \cite{JKS13} and further to Theorem \ref{Forelli-original} if the pencil is standard.

The proof of the theorem of Cho-Kim was established in two steps as follows:

\medskip
\begin{narrower}
	{Step 1.} There is a subpencil on whose underlying set 
	$f$ is holomorphic. 
	
	{Step 2.}  There is a standard subpencil of 
	the pencil found in Step 1.

\end{narrower}
\medskip 
\noindent
First, the analysis of Joo-Kim-Schmalz (\cite{JKS13}, p.1173--1174)  implies that $f$ satisfies the Cauchy-Riemann equations at each point of a Riemann surface of the pencil near the origin. This, together with an application of the Baire category theorem, settles Step 1. If one assumes that Step 2 fails to hold, then a contradiction can be derived from the transversality of each pair of holomorphic curves of the given pencil. Hence we can find a standard subpencil along which the function $f$ is holomorphic. Then we recall the original steps of Forelli and show first that the formal Taylor series $S_f$ is of holomorphic type. Now the conclusion follows from Theorem \ref{CK21} and the generalization of Lemma \ref{Hartogslemma} in \cite{Chirka06}.

\subsection{Vector field viewpoint}

Let $X$ be a diagonalizable vector field of aligned type on $\mathbb{C}^n$. Kim-Poletsky-Schmalz \cite{KPS09} proved that, if a function $f:B^n\to \mathbb{C}$ satisfies: (1) $f\in C^{\infty}(0)$, and (2) $f$ is holomorphic along each integral curve of $X$, then $f$ is holomorphic on $B^n$. For the proof, the authors also followed the original steps of Forelli. Note first that the complex flow map of $X$ is given by
\begin{equation}\label{flow}
\Phi^X(z,t)=(z_1e^{-\lambda_1 t},\dots,z_ne^{-\lambda_n t}).
\end{equation}
If $S_f=\sum C_{IJ}z^I\bar{z}^{J}$ is the formal Taylor series of the function $f$, then Condition (2) implies that the map
\begin{equation}\label{asymptoticexp}
f(\Phi^X(z,t))=\sum \{C_{IJ}z^I\bar{z}^{J}\cdot (\text{exponential term in variable}\;t)\}
\end{equation}
is holomorphic in $t$ for each $z\in B^n$. Here, we used the multi-index notations: 
\begin{align*}
	k=(k_1,\ldots,k_n),~|k|=k_1+\cdots+k_n,~k!=k_1!\cdots k_n!,\,\text{and}\,z^{k}=z_1^{k_1} \cdots z_n^{k_n},
\end{align*}
where $z=(z_1,\dots,z_n)\in \mathbb{C}^n$. Then (\ref{asymptoticexp}) implies that $C_{IJ}=0$ whenever $J\neq 0$; so $S_f$ is of holomorphic type. To prove that $S_f$ is holomorphic on a neighborhood of the origin, the authors of \cite{KPS09} develop a version of the Cauchy estimate along each integral curve of $X$. So they obtain a uniform estimate 
\[
\Big|\sum_{\text{finite}} C_{I0}z^I\Big|<A,~\text{for each}~z\in B^n
\] 
for some $A>0$. Then by applying the Cauchy estimate to the finite polynomial above, one can conclude from the Weierstrass $M$-test that $f\equiv S_f$ is holomorphic on some $B^n(0;r)$. Finally, use the implicit function theorem to `straighten' the complex flow of $X$ locally, and extend $f$ along the flow by Lemma \ref{Hartogslemma}.	

If the given field $X$ is not diagonalizable, then the formal power series expansion of the map $f(\Phi^X(z,\cdot))$ becomes considerably more complicated than (\ref{asymptoticexp}). Therefore it is not clear whether a version of the Cauchy estimate along the leaves of $X$ can be formulated. So the authors of \cite{JKS16} take a different approach: they first use sophisticated mathematical induction on the index $J$ in $S_f=\sum C_{IJ}z^I\bar{z}^J$ to show that $S_f$ is of holomorphic type. Then they apply the Phragm$\acute{\textup{e}}$n-Lindel$\ddot{\textup{o}}$f type argument and the methods in \cite{JKS13} to $f$ and show that $f$ satisfies the Cauchy-Riemann equations along each complex integral curve of $X$. Recall that the set of complex integral curves of $X$ forms a foliation of a punctured open neighborhood $V-\{0\}$ of the origin.  Since $f\in C^1$ on some open neighborhood of the origin, we have $\bar{\partial}f\equiv 0$ on $V-\{0\}$ by shrinking $V$ if necessary. Therefore, $f$ becomes holomorphic on $V$ by Theorem \ref{Hartogs-original} as the origin is a removable singularity of a holomorphic function $f|_{V-\{0\}}$. Arguing as before, the function extends holomorphically on the union of $V$ and the maximal integral curves of $X$ intersecting $V$.

\section{Localization of Kim-Poletsky-Schmalz theorem}

Note that the underlying set of each family of curves considered in the previous section is open. But in Section \ref{Section3}, we only required the underlying set of a suspension of linear discs to be $\textit{nonpluripolar}$ to guarantee $\textbf{Step 2}$. A general nonpluripolar set needs not to have a nonempty interior; indeed, there are nonpluripolar sets even having Lebesgue measure zero (e.g. $\mathbb{R}$ in $\mathbb{C}$). So in the light of the original steps of Forelli, we obtain a new generalization of Theorem \ref{Forelli-original} as soon as we find a nonpluripolar linear suspension with an empty interior on which $\textbf{Step 1}$ is available. 

Motivated by the observation, I proved in \cite{Cho22} that $\textbf{Step 1}$ is valid on a given linear suspension if, and only if, it has a specific leaf $L_v:=\{zv:z\in B^1\}$ $(v\in \bar{F})$ called an $\textit{algebraically}$ $\textit{nonsparse leaf}$. Then I also introduced the potential theoretic notion $\textit{regular leaf}$ and showed that a suspension is nonpluripolar if, and only if, it has a regular leaf. So the Forelli-type theorem holds on a linear suspension having the two kinds of leaves and such a suspension is called a $\textit{Forelli suspension}$. This localization of Theorem \ref{Forelli-original} can be used to construct a nowhere dense Forelli suspension as verified in Example 4.3 of \cite{Cho22}.
 
 Recently, I found that the aforementioned notions of leaves can be generalized to the case of suspensions generated by diagonalizable vector fields of aligned type. Let $X$ be one of such fields on $\mathbb{C}^n$ with eigenvalues $\lambda=(\lambda_1,\dots,\lambda_n)$ and $\Phi^X$ the complex flow map of $X$ in (\ref{flow}). By a $\textit{suspension of} $ $\textit{integral curves of}\;X$,  we mean a pair of the form $(S^X_0(F),\Phi^{X})$. Here, 
 \[
 S^X_0(F):=\{\Phi^{X}(z,t):z\in F,~t\in \mathbb{H}\}
 \]
 for some nonempty subset $F$ of $S^{2n-1}$ and $\mathbb{H}$ is the right open half-plane in $\mathbb{C}.$ We will identify a suspension with its underlying set. Suppose that the given suspension possesses an algebraically nonsparse leaf 
  \[
 L_{z_0}:=\{\Phi^X(z_0,t):t\in \mathbb{H}\},~z_0\in \bar{F}.
 \] 
 Then it turns out that any function $f:B^n\to \mathbb{C}$ satisfying the following two conditions: (1) $f\in C^{\infty}(0)$, and (2) $t\in \mathbb{H}\to f(\Phi^X(z,t))$ is holomorphic for each $z\in F$ has a formal series $S_f$ of holomorphic type. Now $S_f$ can be written as a formal sum of holomorphic $\textit{quasi-homogeneous}$ polynomials of type $\lambda$ (see Definition 1.2 in \cite{Cho22-2}) that converges pointwise on the suspension. If the given suspension has a regular leaf, then the suspension is nonpluripolar so we can use pluripotential methods to establish the uniform convergence of the polynomials. Among several methods available, I chose to follow Siciak's ideas \cite{Siciak82} and develop a new capacity theory which seems more suitable for the analysis of quasi-homogeneous polynomials; see p.4 of \cite{Cho22-2}.
 
 For this purpose, the set of $\textit{quasi-homogeneous}$ plurisubharmonic function on $\mathbb{C}^n$ associated with $X$ was defined in \cite{Cho22-2} as
 	\[
 H_{\lambda}:=\{u\in \textup{PSH}(\mathbb{C}^n): u\geq 0~\text{on}~\mathbb{C}^n, \,u(\Phi^{X}(z,t))=e^{-\textup{Re}\,t}\cdot u(z)\; \forall z\in \mathbb{C}^n, t\in \mathbb{C}\}.
 \]
Then the $\lambda\textit{-projective capacity}$ $\rho_{\lambda}$ and the related extremal function $\Psi_{E,\lambda}$ can be defined as in Definition \ref{defofext}. If $X$ is the complex Euler vector field, i.e., $\lambda=(1,\dots,1)$, then the functions reduce to those in Definition \ref{defofext}.  Using the methods in \cite{Siciak82} and \cite{KPS09}, one can develop a theory of the new functions and obtain the following localization of the theorem of Kim-Poletsky-Schmalz \cite{KPS09}.

\begin{theorem}[Cho \cite{Cho22-2}]\label{main theorem}
	If a suspension $S^X_0(F)$ has a nonsparse leaf and a regular leaf, then it is a Forelli suspension; that is, any function $f:B^n\to \mathbb{C}$ satisfying the following two conditions
	\begin{enumerate}
		\setlength\itemsep{0.1em}
		\item $f\in C^{\infty}(0)$, and
		\item $t\in \mathbb{H}\to f\circ \Phi^{X}(z,t)$ is holomorphic for each $z\in F$
	\end{enumerate}
	is holomorphic on a domain of holomorphy
	\[
	\Omega:=\{z\in \mathbb{C}^n:\Psi^{\ast}_{S^X_0(F),\lambda}(z)<1\}\supset B^n(0;\{\rho_{\lambda}(S^{X}_0(F))\}^{\textup{max}(\lambda)})
	\]
	of the origin. Furthermore, there exists an open neighborhood $U=U(F,X)\subset S^{2n-1}$ of a generator $v_0\in \bar{F}$ of the regular leaf such that $f|_{\Omega}$ extends to a holomorphic function on a domain of holomorphy
	\[
	\hat{\Omega}:=\{z\in \mathbb{C}^n:\Psi_{\Omega \cup S^X_0(U),\lambda}(z)<1\}\supset \Omega \cup S^X_0(U).
	\]
\end{theorem}

 The generalization of Lemma \ref{Hartogslemma} by Shiffman \cite{Shiffman89} plays a crucial role in the extension of $f$ along $S^X_0(U)$; see Section 6 in \cite{Cho22-2}. At this point, a few remarks regarding Theorem \ref{main theorem} should be in order. As each point of a nonempty open subset $U$ of $S^{2n-1}$ generates a leaf that is both nonsparse and regular, $S^X_0(U)$ is always a Forelli suspension. But a Forelli suspension needs not to be generated by an open subset of $S^{2n-1}$ in general; one can construct a set $F\subset S^{2n-1}$ such that $S^X_0(F)$  is a nowhere dense Forelli suspension for any $X$. See Example 7.3 in \cite{Cho22-2}. Note also that we required the given suspension to possess all of the two specific leaves in Theorem \ref{main theorem}; several examples of suspensions in Section 7 of \cite{Cho22-2} indicate that they are indeed all necessary to obtain the desired conclusion.
 
 Finally, we remark that it would be interesting to know whether this type of analysis can be carried out when the given suspension is generated by a nondiagonalizable vector field of aligned type or a family of Riemann surfaces pairwise transversal at the origin.

\vspace{50pt}

Ye-Won Luke Cho (\texttt{ww123hh@pusan.ac.kr}) 

\medskip

Department of Mathematics,

Pusan National University, 

Busan 46241, The Republic of Korea.
\end{document}